\documentclass[twoside, 11pt]{article}
\usepackage[english]{babel}
\usepackage[ansinew]{inputenc}
\usepackage{amsmath,amsthm,amsfonts,amssymb,graphicx}
\usepackage{latexsym}

\newcommand{\ere}{\mathbb{R}} 
\newcommand{\fle}{\rightarrow}

 \newtheorem{theorem}{Theorem}[section]
 \newtheorem{cor}[theorem]{Corollary}
 \newtheorem{lemma}[theorem]{Lemma}
 
 \theoremstyle{definition}
 
 \theoremstyle{remark}
 \newtheorem{remark}[theorem]{Remark}
 \theoremstyle{eg}

 \theoremstyle{fact}
 
 \theoremstyle{remark}
  \newtheorem*{claim*}{Claim}
  \theoremstyle{remark}
  
\numberwithin{equation}{section}

\newcommand{\bthm}{\begin{theorem}}
\newcommand{\ethm}{\end{theorem}}
\newcommand{\bprop}{\begin{proposition}}
\newcommand{\eprop}{\end{proposition}}
\newcommand{\bcor}{\begin{corollary}}
\newcommand{\ecor}{\end{corollary}}
\newcommand{\blem}{\begin{lemma}}
\newcommand{\elem}{\end{lemma}}

\newcommand{\ba}{\begin{array}}
\newcommand{\ea}{\end{array}}
\newcommand{\be}{\begin{enumerate}}
\newcommand{\ee}{\end{enumerate}}
\newcommand{\beq}{\begin{equation}}
\newcommand{\eeq}{\end{equation}}
\newcommand {\bea} {\begin{eqnarray}}
\newcommand {\eea} {\end {eqnarray}}
\newcommand {\bua} {\begin{eqnarray*}}
\newcommand {\eua} {\end {eqnarray*}}

\newcommand{\ds}{\displaystyle}
\newcommand{\se}{\subseteq}

\newcommand{\N}{{\mathbb N}}
\newcommand{\Z}{{\mathbb Z}}
\newcommand{\vp}{\varphi}
\newcommand{\eps}{\varepsilon}

\newcommand{\zgeqk}{\{k,k+1,\ldots\}}

\title{\bf Alternative iterative methods for nonexpansive mappings, rates of convergence and applications}

\date{}
\begin{document}
\maketitle

\begin{center}
{\bf Vittorio Colao$^{a,}$\footnote{Supported by Regione Calabria POR/FSE 2007-2013 and Universit\'a della Calabria. E-mail address: colao@mat.unical.it}, Lauren\c tiu Leu\c
stean$^{b,c,}$\footnote{Supported by the Oberwolfach Leibniz Fellow
Programme (OWLF) at the Oberwolfach Mathematics Institute and by the
German Science Foundation (DFG Project KO 1737/5-1). E-mail address:
leustean@mathematik.tu-darmstadt.de}, Genaro L\' opez$^{d,}$\footnote{Supported by DGES, Grant MTM2006-13997-C02-01 and Junta de
Andalucía, Grant FQM-127 and partially by  the Oberwolfach
Mathematics Institute under the OWLF Programme. E-mail address:
glopez@us.es},\\ Victoria Martín-Márquez$^{d,}$\footnote{Supported
by Junta de Andalucía, Grant FQM-127 and Ministerio de Ciencia e
Innovación, Grant AP2005-1018. E-mail address: victoriam@us.es}}
\end{center}
$^a${\small Dipartimento di Matematica, Universit\'a della Calabria,
87036 Arcavacata di Rende (CS), Italy}\\
$^b${\small Department of Mathematics, Technische Universität Darmstadt,
 Schlossgartenstrasse 7, 64289 Darmstadt, Germany}\\
$^c$ {\small Institute of Mathematics "Simion Stoilow'' of the
Romanian Academy, Calea Grivi\c tei 21, 010702 Bucharest, Romania}\\
$^d${\small Departamento de Análisis Matemático, Universidad de
Sevilla, Apdo. 1160, 41080-Sevilla, Spain}

\vskip.5in

\begin{abstract}
Alternative iterative methods for a nonexpansive mapping in a Banach
space are proposed and proved to be convergent to a common solution
to a fixed point problem and a variational inequality. We give rates
of asymptotic regularity for such iterations using proof-theoretic
techniques. Some applications of the convergence results are
presented.


\end{abstract}

\section{Introduction}
Many problems arising in different areas of mathematics such as
optimization, variational analysis and game theory, can be
formulated as the fixed point problem:
\begin{equation}\label{FPP}
\textrm{ find } x\in X \textrm{ such that } x=Tx,
\end{equation}
where $T$ is a nonexpansive mapping defined on a metric space $X$,
i.e., $T$ satisfies the property  $d(Tx,Ty) \leq d(x,y) $, for all
$x,y\in X$.

For instance, let $A:C\to H$ be a nonlinear operator where $C\subset
H$ is a closed convex subset of a Hilbert space. The variational
inequality problem associated to $A$, VIP($A,C$), is formulated as
finding a point $x^* \in C$ such that
\begin{equation}\label{VIP-intro}
 \langle Ax^*  , x-x^* \rangle\geq 0 \,\, \,\, \forall x\in C .
\end{equation}
It is well-known that the VIP($A,C$) is equivalent to the problem
of finding the fixed point
\begin{equation}\label{FPP.intro}
 x^* = P_C ( x^* - \lambda Ax^*),
\end{equation}
where $ \lambda >0$ and $P_C$ is the metric projection onto $C$,
which is a nonexpansive mapping in this case. Besides, if $ f:C\to
\ere$ is a differentiable convex function and we denote by $A$ the
gradient operator of $f$, then \eqref{VIP-intro} is the optimality
condition for the minimization problem
\begin{equation}\label{MP-intro}
\min_{x\in C} f(x).
\end{equation}

Bearing in mind that the iterative methods for approximating a fixed
point of a nonexpansive mapping can be applied to find a solution to
a variational inequality, zeros of an accretive operator and a
minimizer of a convex function, in the recent years the study of the
convergence of those methods has received a great deal of attention.
Basically two types of iterative algorithms have been investigated:
Mann algorithm and Halpern algorithm.

In the following, let $X$ be a real Banach space, $C\subset X$ a
closed convex subset and $T:C\to C$ a nonexpansive mapping with
fixed point set $F=\{x\in C:x=Tx\}\neq\emptyset$.

Mann algorithm generates a sequence according to the following
recursive manner:
\begin{equation}\label{MI}
x_{n+1}=(1-t_n)x_n+t_nTx_n,\quad n\ge 0,
\end{equation}
where the initial guess $x_0\in C$ and $\{t_n\}$ is a sequence in
$(0,1)$.

Halpern algorithm generates a sequence via the recursive formula:
\begin{equation}\label{HI}
 x_{n+1} =\alpha_n u +(1-\alpha_n ) T x_n ,\quad  n\geq 0
\end{equation}
where $x_0, u\in C$ are arbitrary and the sequence
$\{\alpha_n\}\subset(0,1)$.

Whenever a fixed point of the mapping $T$ exists, while Halpern
algorithm strongly converges, we just get weak convergence for Mann
algorithm, as was established in \cite{GL75} thanks to a
counterexample. The references \cite{M53,K55,I76,R79,GKKR01,X03} can
be consulted for convergence results of Mann algorithm. Some
modifications have been proposed in \cite{NT03,KX05} to get strong
convergence. As for Halpern algorithm, see
\cite{H67,L77,W92,R83,X02,CC06,S07,LMX} and references therein for
studies dedicated to its convergence.

Another iterative approach to solving the problem \eqref{FPP}
which may have multiple solutions, is to replace it by a family of
perturbed problems admitting a unique solution, and then to get a
particular original solution as the limit of these perturbed
solutions as the perturbation vanishes. For example, Browder
\cite{B66,B67} proved that if the underlying space $H$ is Hilbert,
then, given $u\in H$ and $t\in(0,1)$, the approximating curve
$\{x_t\}$ defined by
\begin{equation}\label{BI}
x_t=tu+(1-t)Tx_t
\end{equation}
strongly converges, as $t\to 0$, to the fixed point of $T$ closest
to $u$ from $F$. Browder's result has been generalized and
extended to a more general class of Banach spaces
\cite{R80,ST97,OPX06}. Combettes and Hirstoaga \cite{CH06}
introduced a new type of approximating curve for fixed point
problems in the setting of a Hilbert space. This curve defined by
the implicit formula
\begin{equation}\label{CHI}
x_t = T (tu+(1-t)x_t),
\end{equation}
was proved to converge to the best approximation to $u$ from $F$. In
\cite{X09} Xu studied the behavior of $\{x_t\}$ defined by
\eqref{CHI} in the setting of a Banach space $X$ and discretized
this regularization method studying the strong convergence of the
explicit algorithm
\begin{equation}\label{XI}
x_{n+1} = T (\alpha_n u +(1-\alpha_n)x_n),
\end{equation}
where $\{\alpha_n\}\subset(0,1)$. Moreover, he proved that the
convergence point is the image of $u$ under the unique sunny
nonexpansive retraction $Q$ from $X$ to $F$ (see, for instance,
\cite{R73, R80}).

On the other hand, Moudafi in \cite{M00} introduced the viscosity
approximation method for nonexpansive mappings, which generalizes
Browder's \eqref{BI} and Halpern \eqref{HI} iterations, by using a
contraction $\Phi$ instead of an arbitrary point $u$. The
convergence of the implicit and explicit algorithms has been the
subject of many papers because under suitable conditions these
iterations strongly converge to the unique solution $q\in F$ to the
variational inequality
\begin{equation}\label{vi}
\langle(I-\Phi)q, J(x-q)\rangle \geq 0 \,\, \,\, \forall x\in F,
\end{equation}
where $J$ is a duality mapping, i.e., $q$ is the unique fixed
point of the contraction $Q\circ \Phi$. This fact allows us to
apply this method to convex optimization, linear programming and
monotone inclusions. See \cite{X04,SC06,S07b} and references
therein for convergence results regarding viscosity approximation
methods.

In this paper, we analyze the behavior of a new approximating
curve in the setting of Banach spaces, which constitutes a hybrid
method of the ones presented by Combettes and Hirstoaga
\eqref{CHI} and Moudafi. This curve is defined by
\begin{equation}\label{Imp Iter}
x_t = T (t\Phi(x_t)+(1-t)x_t),
\end{equation}
for some contraction $\Phi$, that is, for any $t\in(0,1)$ $x_t$ is
the unique fixed point of the contraction $T_t=T(t\Phi+(1-t)I)$.
The discretized iteration
\begin{equation}\label{Exp Iter}
x_{n+1} = T(\alpha_n\Phi(x_n)+(1-\alpha_n)x_n),
\end{equation}
is also considered and studied under suitable conditions on the
sequence $\{\alpha_n\}\subset(0,1)$. From this explicit algorithm
we obtain the so-called hybrid steepest descent method
\begin{equation}\label{hybrid}
x_{n+1}=Tx_n-\alpha_n g(Tx_n).
\end{equation}
This iterative method was suggested by Yamada \cite{Y01} as an
extension of viscosity approximation methods for solving the
variational inequality VIP($g,F$) \eqref{VIP-intro} in the case
when $g$ is strongly monotone and Lipschitz continuous, and $F$ is
the fixed point set of a mapping $T$ which belongs to a subclass
of the quasi-nonexpansive mappings (also see \cite{Y04,M08}). We
will get the convergence of the algorithm \eqref{hybrid} for a
nonexpansive mapping $T$, just requiring $I-\mu g$ to be a
contraction for some $\mu>0$, which it is satisfied in the
particular case when $g$ is strongly monotone and Lipschitz
continuous.

Asymptotic regularity is a very important concept in metric fixed
point theory. It was already implicit in \cite{K55,S57,E70}, but it
was formally introduced by Browder and Petryshyn in \cite{BP66}. In
our setting, the mapping $T$ is called asymptotically regular if for
all $x\in C$
$$\lim_{n\to\infty} \Vert x_n -Tx_n \Vert =0.$$
Effective rates of asymptotic regularity for both Mann and Halpern
iterations have been obtained (see \cite{K01,K03,KL03,L07a,L07b}) by
applying methods of proof mining. By ``proof mining'' we mean the
logical analysis, using proof-theoretic tools, of mathematical
proofs with the aim of extracting relevant information hidden in the
proofs. This new information can be both of quantitative nature,
such as algorithms and effective bounds, as well as of qualitative
nature, such as uniformities in the bounds or weakening the
premises. Thus, even if one is not particularly interested in the
numerical details of the bounds themselves, in many cases such
explicit bounds immediately show the independence of the quantity in
question from certain input data. A comprehensive reference for
proof mining is Kohlenbach's  book \cite{K08}. One of the aims of
this paper is to give effective rates of asymptotic regularity for
the algorithm \eqref{Exp Iter} of nonexpansive mappings in the
framework of normed spaces.

The organization of the paper is as follows. In section 2 we
introduce some preliminary results and present a technical lemma
with regards to the behavior of the sequence defined by the
algorithm \eqref{Exp Iter}, which will be useful for the proof of
the convergence of the sequence and the evaluation of the rate of
asymptotically regularity in following sections. Section 3 contains
the main results about the strong convergence of both implicit
\eqref{Imp Iter} and explicit \eqref{Exp Iter} algorithms to the
unique solution to the variational inequality \eqref{vi} in the
setting of uniformly smooth Banach spaces, and also in the framework
of reflexive Banach spaces with weakly continuous normalized duality
mapping in the case of the implicit iteration. Section 4 is devoted
to the rate of asymptotic regularity for the iterations \eqref{Exp
Iter} of nonexpansive mappings in normed spaces. Finally, in section
5 we give examples of how to apply the main results of section 3 to
find a solution to a variational inequality or a zero of an
accretive operator.

\section{Preliminaries}

Let $X$ be a real Banach space with norm $\|\cdot\|$ and dual
space $X^*$. For any $x\in X$ and $x^*\in X^*$ we denote
$x^*(x)=\langle x, x^*\rangle$. Given a nonempty closed convex
subset $C\subset X$, $\Phi:C\fle C$ will be a a $\rho$-contraction and
$T:C\rightarrow C$ a nonexpansive self-mapping with nonempty fixed
point set $F:=\{x\in C:Tx=x\}$.

We include some brief knowledge about geometry of Banach spaces
which can be found in more details in \cite{C90}. The
\emph{normalized duality mapping} $J:X\fle 2^{X^*}$ is defined by
\begin{equation}
J(x)=\{x^*\in X^{^*} : \langle x,x^* \rangle = \|x\|^2 =
\|x^*\|^2\}.
\end{equation}
It is known that
$$J(x) = \partial(\|x\|),$$
where $\partial (\|x\|)$ is the subdifferential of $ \|\cdot\|$ at
$x$ in the sense of convex analysis. Thus, for any $x,y \in X$, we
have the subdifferential inequality
\begin{equation}\label{subdiff ineqN}
\|x+y\|^2 \leq \|x\|^2+ 2\langle y, j(x+y)\rangle, \quad j(x+y)\in
J(x+y).
\end{equation}
A Banach space $X$ is said to be \emph{smooth} if
\begin{equation}\label{lim}
\lim_{t\fle0} \frac{\|x+ty\|-\|x\|}{t}
\end{equation}
exists for each $x, y \in S_X$, where $S_X$ is the unit sphere of
$X$, i.e., $S_X = \{v\in X : \|v\| = 1\}$. When this is the case,
the norm of $X$ is said to be \emph{Gâteaux differentiable}. If for
each $y\in X$ the limit \eqref{lim} is uniformly attained for $x\in
X$, we say that the norm of $X$ is \emph{uniformly Gâteaux
differentiable}, and we say that $X$ is \emph{uniformly smooth} if
the limit \eqref{lim} is attained uniformly for any $x, y\in S_X$.

It is known that a Banach space $X$ is smooth if and only if the
duality mapping $J$ is single-valued, and that $X$ is uniformly
smooth if and only if the duality mapping $J$ is single-valued and
norm-to-norm uniformly continuous on bounded sets of $X$.
Moreover, if $X$ has a uniformly Gâteaux differentiable norm then
$J$ is norm-to-weak$^*$ uniformly continuous on bounded sets of
$X$.

Following Browder \cite{B67} we say that the duality mapping $J$
is \emph{weakly sequentially continuous} (or simply \emph{weakly
continuous}) if $J$ is single-valued and weak-to-weak$^*$
sequentially continuous; i.e., if $x_n \rightharpoonup x$ in $X$,
then $J(x_n)\rightharpoonup^*J(x)$ in $X^*$. A Banach space with
weakly continuous duality mapping is known (see \cite{LX94}) to
satisfy \emph{Opial's property} (i.e., whenever $x_n
\rightharpoonup x$ and $y\neq x$, we have $\overline{\lim} \|x_n
-x\| < \overline{\lim} \|x_n - y\|$), and this fact implies (see
\cite{GK90}) that $X$ satisfies the \emph{Demiclosedness
principle}: if $C$ is a closed convex subset of $X$ and $T$ is a
nonexpansive self-mapping, then $x_n\rightharpoonup x$ and
\mbox{$(I-T)x_n \rightarrow y$} imply that $(I-T)x=y$ .

Consider a subset $D\subset C$ and a mapping $Q:C\fle D$. We say
that $Q$ is a \emph{retraction} provided $Qx=x$ for any $x\in D$.
The retraction $Q$ is said to be \emph{sunny} if it satisfies the
property: $Q(x+t(x-Qx))=Qx$ whenever $x+t(x-Qx)\in C$, where $x\in
C$ and $t\geq 0$.

\begin{lemma} \cite{Br73, R73, GK84} \label{sunny}
Let $X$ be a smooth Banach space and $D \subset C$ be nonempty
closed convex subsets of $X$. Given a retraction $Q: C \rightarrow
D$, the following three statements are equivalent:
\begin{itemize}
\item [(a)] $Q$ is sunny and nonexpansive.

\item [(b)] $\|Qx-Qy\|^2\leq\langle x-y,J(Qx-Qy)\rangle$ for all
$x,y\in C$.

\item [(c)] $\langle x-Qx,J(y-Qx)\rangle\leq0$ for all $x\in C$
and $y\in D$.
\end{itemize}
Consequently, there is at most one sunny nonexpansive retraction
from $C$ onto $D$.
\end{lemma}

In some circumstances, we can construct the sunny nonexpansive
retraction. For the nonexpansive mapping $T$ with fixed point set
$F$, an arbitrary $u\in C$ and $t\in (0,1)$, let $z_t$ be the unique
fixed point of the contraction $z\mapsto tu+(1-t)Tz$ for $z\in C$;
that is, $z_t$ is the unique solution in $C$ to the fixed point
equation:
\begin{equation}\label{eq:2:fp}
z_t=tu+(1-t)Tz_t.
\end{equation}
It is natural to study the behavior of the net $\{z_t\}$ as $t\to
0^+$. However, it is unclear if the strong $\lim_{t\to 0^+}z_t$
always exists in a general Banach space. The answer is yes in some
classes of smooth Banach spaces and then the limit defines the
sunny nonexpansive retraction from $C$ onto $F$. Those Banach
spaces where the net $\{z_t\}$ strongly converges are said to have
\emph{Reich's property} since Reich was the first to show that all
uniformly smooth Banach spaces have this property.

\begin{theorem} \cite{R80,OPX06} \label{unique sunny1}
Let $X$ be either a uniformly smooth Banach space or a reflexive
Banach space with a weakly continuous duality mapping, $C$ be a
nonempty closed convex subset of $X$, and  $T: C\to C$ be a
nonexpansive mapping with $F\neq\emptyset$. Then the net $\{z_t\}$
strongly converges as $t\to 0^+$ to a fixed point of $T$;
moreover, the limit
\begin{equation}\label{eq:2:retrac}
Q(u):=\lim_{t\rightarrow0^+} z_t
\end{equation}
defines the unique sunny nonexpansive retraction from $C$ onto
$F$.
\end{theorem}

In \cite{R80b} Reich proved the following two lemmas which will be
needed for the convergence results in section 3.
\begin{lemma}
\label{lem:reich1}Let $\{x_{n}\}$ be a bounded sequence contained
in a separable subset $D$ of a Banach space $X.$ Then there is a
subsequence $\{x_{n_{k}}\}$ of $\{x_{n}\}$ such that
\[
\lim_{k}\|x_{n_{k}}-y\|\] exists for all $y\in D.$
\end{lemma}

\begin{lemma}
\label{lem:reich2}Let $D$ be a closed convex subset a real Banach
space $X$ with a uniformly Gâteaux differentiable norm, and let
$\{x_{n}\}$ be a sequence in $D$ such that \[
f(y):=\lim_{n}\|x_{n}-y\|\] exists for all $y\in D$. If $f$
attains its minimum over $D$ at $u$, then \[ \limsup_{n}\langle
y-u,j(x_{n}-u)\rangle\leq0\]
 for all $y\in D.$
\end{lemma}

The following lemma collects some properties of the iteration
\eqref{Exp Iter}, useful both for proving the convergence of the
iteration and for computing the rate of asymptotic regularity.

\begin{lemma}\label{lemma-useful}
Let $X$ be a normed space and $\{x_n\}$ be the sequence defined by
the explicit algorithm \eqref{Exp Iter}.
\begin{itemize}
\item[(1)] For all $n\ge 0$,
\begin{eqnarray}
\|\Phi(x_n)-x_n\| &\le & (1+\rho)\|x_n-x_0\|+\|\Phi(x_0)-x_0\|\label{ineq-Phi-xn},\\
\|x_n-Tx_n\| &\le & \|x_{n+1}-x_n\| +
\alpha_n\|\Phi(x_n)-x_n\|.\label{ineq-xn-Txn}
\end{eqnarray}
\item[(2)] For all $n\ge 1$,
\beq \ba{lll}
\|x_{n+1}-x_n\| &\le & (1-(1-\rho)\alpha_n)\|x_n-x_{n-1}\|\\[0.1cm]
&&\,\, +\,|\alpha_n-\alpha_{n-1}|\cdot\|\Phi(x_{n-1})-x_{n-1}\|. \ea
\label{ineq-xn+1-xn} \eeq
\item[(3)] If $T$ has fixed points, then $\{x_n\}$ is bounded for every $x_0\in C$.
\end{itemize}
\end{lemma}

\begin{proof}$\,$
\begin{itemize}
\item[(1)] Let $n\ge 0$.
\bua
\|\Phi(x_n)-x_n\| &\le& \|\Phi(x_n)-\Phi(x_0)\|+\|\Phi(x_0)-x_0\|+\|x_0-x_n\|\\
&\le& (1+\rho)\|x_n-x_0\|+\|\Phi(x_0)-x_0\|.\\
\|x_n-Tx_n\| &\le & \|x_{n+1}-x_n\| +\|x_{n+1}-Tx_n\|\\
&=& \|x_{n+1}-x_n\| +\|T(\alpha_{n}\Phi(x_n)+(1-\alpha_n)x_n)-Tx_n\|\\
&\le & \|x_{n+1}-x_n\| + \|\alpha_{n}\Phi(x_n)+(1-\alpha_n)x_n-x_n\|\\
&=&\|x_{n+1}-x_n\| +\alpha_n\|\Phi(x_n)-x_n\|. \eua
\item[(2)] Let $n\ge 1$.
\bua
\|x_{n+1}-x_n\| &=&\|T(\alpha_{n}\Phi(x_n)+(1-\alpha_n)x_n)\\
&& \,\,-\, T(\alpha_{n-1}\Phi(x_{n-1})+(1-\alpha_{n-1})x_{n-1})\|\\
&\le &\|\alpha_{n}\Phi(x_n)+(1-\alpha_n)x_n-\alpha_{n-1}\Phi(x_{n-1})\\
&& \,\,-\,(1-\alpha_{n-1})x_{n-1}\|\\
&=& \|\alpha_n(\Phi(x_n)-\Phi(x_{n-1}))+(1-\alpha_n)(x_n-x_{n-1})\\
&&\,\,+\,(\alpha_n-\alpha_{n-1})(\Phi(x_{n-1})-x_{n-1})\|\\
&\le & \alpha_n\rho\|x_n-x_{n-1}\|+(1-\alpha_n)\|x_n-x_{n-1}\|\\
&&\,\,+\,|\alpha_n-\alpha_{n-1}|\cdot\|\Phi(x_{n-1})-x_{n-1}\|\\
&=& (1-(1-\rho)\alpha_n)\|x_n-x_{n-1}\|\\
&& \,\,+\,|\alpha_n-\alpha_{n-1}|\cdot\|\Phi(x_{n-1})-x_{n-1}\|.
\eua
\item[(3)] Let $p$ be a fixed point of $T$.
\bua
\|x_{n+1}-p\| & = & \|T(\alpha_{n}\Phi(x_n)+(1-\alpha_n)x_n)-Tp\|\\
&\le &\|\alpha_{n}\Phi(x_n)+(1-\alpha_n)x_n-p\|\\
&=& \|\alpha_{n}(\Phi(x_n)-\Phi(p))+(1-\alpha_n)(x_n-p)+\alpha_{n}(\Phi(p)-p)\|\\
&\le & \alpha_n\rho\|x_{n}-p\|+(1-\alpha_n)\|x_{n}-p\|+\alpha_{n}\|\Phi(p)-p\|\\
&=& (1-(1-\rho)\alpha_n)\|x_{n}-p\|+(1-\rho)\alpha_n\frac{\|\Phi(p)-p\|}{1-\rho}\\
&\le & \max\left\{\|x_{n}-p\|, \frac{\|\Phi(p)-p\|}{1-\rho}\right\}.
\eua By induction, we obtain that for all $n\ge 0$,
\[\|x_{n}-p\|\le \max\left\{\|x_0-p\|, \frac{\|\Phi(p)-p\|}{1-\rho}\right\},\]
thus $\{x_n\}$ is bounded. \end{itemize}
\end{proof}

Let us recall some terminology that is used for expressing the
quantitative results in Section \ref{s-rate-as-reg}. We denote by
$\Z_+$ the set of nonnegative integers. Let $k\in\Z_+$ and
$\{a_n\}_{n\geq k}$ be a sequence of nonnegative real numbers. If
$\{a_n\}$ is convergent, then a function
$\omega:(0,\infty)\to\zgeqk$ is called a {\em Cauchy modulus} of
$\{a_n\}$ if for all $\eps>0$,
\begin{equation}\label{def-mod-Cauchy}
|a_{\omega(\eps)+n}-a_{\omega(\eps)}| < \eps,\,\, \,\forall
n\in\Z_+.
\end{equation}
If $\ds\lim _{n\to\infty}a_n=a$, then a function
$\omega:(0,\infty)\to\zgeqk$ is called a {\em rate of convergence}
of $\{a_n\}$ if for any $\eps>0$
\begin{equation}\label{def-rate-conv}
|a_{n}-a| < \eps,\,\,\,\forall n\geq\omega(\eps).
\end{equation}
If the series $\ds \sum_{n=k}^\infty a_n$ is divergent, then a
function $\omega:\Z_+\to\zgeqk$ is called a {\em rate of divergence}
of the series if $\ds \sum_{i=k}^{\omega(n)}a_i \geq n$ for all $
n\in\Z_+$. If the series $\ds \sum_{n=k}^\infty a_n$ converges, then
by a {\em Cauchy modulus} of the series we mean a Cauchy modulus of
the sequence of partial sums $\ds \{s_n\}_{n\ge k}, \,
s_n=\ds\sum_{i=k}^n a_i$.

\begin{lemma}\cite{X02}\label{lem:xu1} Assume $\{a_n\}$ is a sequence of nonnegative real
number such that\[
a_{n+1}\leq(1-\gamma_{n})a_{n}+\gamma_{n}b_{n}+\epsilon_{n},\qquad
n\geq0,\] where $\{\gamma_{n}\}$ and $\{\epsilon_{n}\}$ are
sequences in $(0,1)$ and $\{b_{n}\}$ is a sequence in $\mathbb{R}$
such that $\ds\sum_{n=1}^\infty\gamma_{n}=\infty$, $\ds\sum_{n=1}^\infty\epsilon_{n}<\infty$
and either $\limsup_{n}b_{n}\leq0$ or $\ds \sum_{n=1}^\infty\gamma_{n}|b_{n}|<\infty$.

Then $\ds\lim_{n\to\infty} a_n =0$
\end{lemma}

The following lemma is a quantitative version of \cite[Lemma 2]{L95} and has been proved in \cite{L07b}.

\begin{lemma}\label{quant-Liu-L}\cite[Lemma 9]{L07b}\\
Let $\{\lambda_n\}_{n\ge 1}$ be a sequence in $(0,1)$ and $\{a_n\}_{n\geq 1},\{b_n\}_{n\geq 1}$
be sequences of nonnegative real numbers such that
\[
a_{n+1}\leq (1-\lambda_{n+1}) a_n + b_n \quad  \text{for  all } n\in\N.
\]
Assume that $\ds \sum_{n=1}^\infty \lambda_n$ is divergent, $\,\ds \sum_{n=1}^\infty b_n$ is
convergent and let $\delta:\Z_+\to\N$ be a rate of divergence of $\ds \sum_{n=1}^\infty \lambda_n$
and $\gamma:(0,\infty)\to\N$ be a Cauchy modulus of $\,\ds \sum_{n=1}^\infty b_n$.

Then $\ds\lim_{n\to\infty} a_n =0$ and moreover for all $\eps\in(0,2)$
\beq
 a_n < \eps, \quad \forall n\ge
h(\gamma,\delta,D,\eps),\label{conclusion-quant-liu}
\eeq
where $D>0$ is an upper bound on $\{a_n\}$ and
\[
\ds h(\gamma,\delta,D,\eps)=\delta\left(\gamma\left(\frac\eps 2\right)+1+\left\lceil\ln
\left(\frac{2D}\eps\right)\right\rceil\right).
\]
\end{lemma}


\section{Convergence of the algorithms}\label{section3}

In this section we prove the convergence of the implicit \eqref{Imp
Iter}, explicit \eqref{Exp Iter} and hybrid steepest descent
\eqref{hybrid} algorithms in the setting of Banach spaces, which
generalize previous results by Combettes and Hirstoaga \cite{CH06},
Xu \cite{X09}, Yamada \cite{Y01}, and Xu and Kim \cite{XK03}.

\begin{theorem}\label{thm:implicit}
Let $X$ be either a reflexive space with weakly continuous
normalized duality mapping $J$ or a uniformly smooth Banach space,
$C$ a nonempty closed convex subset of $X$, $T:C\to C$ a
nonexpansive mapping with fixed point set $F\neq\emptyset$ and
$\Phi:C\to C$ a $\rho$-contraction. Then the approximating curve
$\{x_t\}\subset C$ defined by
\begin{equation}\label{Implicit Iter}
x_{t}=T(t\Phi(x_t)+(1-t)x_t)
\end{equation}
strongly converges, as $t\fle0$, to the unique solution $q\in F$ to
the inequality
\begin{equation}\label{inequality Th.}
\langle(\Phi-I)q,J(x-q)\rangle\leq0,\,\,\, \forall x\in F.
\end{equation}
\end{theorem}

\begin{proof}
We observe that we may assume that $C$ is separable. To see this,
consider the set $K$ defined as\[
\begin{array}{l}
K_{0}:=\{q\},\\
K_{n+1}:=co(K_{n}\cup T(K_{n})\cup\Phi(K_{n})),\\
K:=\overline{\bigcup_{n}K_{n}}.\end{array}\] Then $K\subset C$ is
nonempty, convex, closed, and separable. Moreover $K$ maps into
itself under $T$, $\Phi$ and, therefore, $T_t=T(t\Phi+(1-t)I)$.
Then $\{x_t\}\subset K$ and we may replace $C$ with $K$.

We will prove that $\{x_t\}$ converges, as $t\fle 0$, to the point
$q\in F$ which is the unique solution to the inequality
\eqref{inequality Th.}.

The sequence $\{x_t\}$ is bounded. Indeed, given $p\in F$,
\[
\begin{array}{ll}
\|x_t-p\| & =\|T(t\Phi(x_t)+(1-t)x_t)-Tp\|\\
 & \le\|t(\Phi(x_t)-\Phi(p))+(1-t)(x_t-p)+t(\Phi(p)-p)\|\\
 & \leq(t\rho+(1-t))\|x_t-p\|+t\|\Phi(p)-p\|.\end{array}\]
Then, for any $t\in (0,1)$,
$$\|x_t-p\|\le\dfrac{1}{1-\rho}\|\Phi(p)-p\|.$$
Take an arbitrary sequence $\{t_n\}\subset (0,1)$ such that
$t_n\fle 0$, as $n\fle0$, and denote $x_n=x_{t_n}$ for any
$n\geq0$. Let
$\Gamma:=\limsup_{n\fle\infty}\langle\Phi(q)-q,J(x_n-q)\rangle$
and $\{x_{n_{k}}\}$ be a subsequence of $\{x_n\}$ such that
$$\lim_{k\fle\infty}\langle\Phi(q)-q,J(x_{n_{k}}-q)\rangle=\Gamma.$$
Since $\{x_{n_{k}}\}$ is bounded, by Lemma \ref{lem:reich1}, there
exists a subsequence, which also will be denoted by
$\{x_{n_{k}}\}$ for the sake of simplicity, satisfying that
$$f(x):=\lim_{k}\|x_{n_{k}}-x\|$$
exists for all $x\in C.$

We define the set
$$A:=\{z\in C:f(z)=\min_{x\in C}f(x)\}$$
and note that $A$ is a nonempty bounded, closed and convex set
since $f$ is a continuous convex function and
$\lim_{\|x\|\to\infty}f(x)=\infty.$ Moreover,\[
\|x_{n_{k}}-Tz\|\leq
t_{n_{k}}\|\Phi(x_{n_{k}})-x_{n_{k}}\|+\|x_{n_{k}}-z\|,\] for any
$z\in C$. This is enough to prove that $T$ maps $A$ into itself.

Since $A$ is a nonempty bounded, closed and convex subset of
either a reflexive space with weakly continuous normalized duality
mapping or a uniformly smooth Banach space, it has the fixed point
property for nonexpansive mappings (see \cite{GK90}), that is
$F\cap A\neq\emptyset.$

If $X$ is reflexive with weakly continuous normalized duality
mapping, we can assume that $\{x_{n_{k}}\}$ has been chosen to be
weakly convergent to a point $\tilde q.$ Since $X$ satisfies
Opial's property, we have $A=\{\tilde q\}$. Then, since $\tilde q \in F$, we obtain
that
$$\Gamma= \langle\Phi(q)-q,J(\tilde q-q)\rangle\leq 0.$$

If $X$ is uniformly smooth, let $\tilde{q}\in F\cap A.$ Then
$\tilde{q}$ minimize $f$ over $C$ and, since the norm is uniformly
Gâteaux differentiable, by Lemma \ref{lem:reich2},
\begin{equation}
\limsup_{k}\langle
x-\tilde{q},J(x_{n_{k}}-\tilde{q})\rangle\leq0\label{eq:minimo}
\end{equation}
holds, for all $x\in C$, and in particular for
$x=\Phi(\tilde{q}).$

We shall show that $\{x_{n_{k}}\}$ strongly converges to
$\tilde{q}.$ Set
$\delta_{k}:=\langle\Phi(\tilde{q})-\tilde{q},J(t_{n_{k}}(\Phi(x_{n_{k}})-x_{n_{k}})+
(x_{n_{k}}-\tilde{q}))-J(x_{n_{k}}-\tilde{q})\rangle.$
Since $J$ is norm-to-weak$^*$ uniformly continuous, we see that
$\lim_{k}\delta_{k}=0.$ Moreover,\begin{equation}
\begin{array}{ll}
\|x_{n_{k}}-\tilde{q}\|^{2} & \leq\|t_{n_{k}}(\Phi(x_{n_{k}})-\Phi(\tilde{q}))+(1-t_{n_{k}})(x_{n_{k}}-\tilde{q})\\
 & \quad+t_{n_{k}}(\Phi(\tilde{q})-\tilde{q})\|^{2}\\
 & \leq\|t_{n_{k}}(\Phi(x_{n_{k}})-\Phi(\tilde{q}))+(1-t_{n_{k}})(x_{n_{k}}-\tilde{q})\|^{2}\\
 & \quad+2t_{n_{k}}\delta_{k}\\
 & \quad+2t_{n_{k}}\langle\Phi(\tilde{q})-\tilde{q},J(x_{n_{k}}-\tilde{q})\rangle\\
 & \leq(1-(1-\rho)t_{n_{k}})\|x_{n_{k}}-\tilde{q}\|^{2}\\
 & \quad+2t_{n_{k}}\delta_{k}\\
 & \quad+2t_{n_{k}}\langle\Phi(\tilde{q})-\tilde{q},J(x_{n_{k}}-\tilde{q})\rangle.
 \end{array}\label{eq:conv_impli}\end{equation}
From \eqref{eq:conv_impli} and by \eqref{eq:minimo}, we obtain\[
\begin{array}{ll}
\lim_{k}\|x_{n_{k}}-\tilde{q}\|^{2} & \leq\limsup_{k}\dfrac{2}{1-\rho}\Bigl(\delta_{k}+
2\langle\Phi(\tilde{q})-\tilde{q},J(x_{n_{k}}-\tilde{q})\rangle\Bigr)\\
 & \leq0.\end{array}\]
That is $\ds \lim_kx_{n_{k}}=\tilde{q}.$ Since $\tilde{q}$ is a
fixed point of $T$, we also have\[
\Gamma=\lim_{k}\langle\Phi(q)-q,J(x_{n_{k}}-q)\rangle=\langle\Phi(q)-q,J(\tilde{q}-q)\rangle\leq0.\]
By applying \eqref{eq:conv_impli} to $\{x_{n}\}$ and $q$, since
$\Gamma\leq 0$ in both cases, we
obtain\[ \lim_{n}x_{n}=q\] as required.
\end{proof}

\begin{remark}
It is easily seen that the conclusion of Theorem
\ref{thm:implicit} remains true if the uniform smoothness
assumption of $X$ is replaced with the following two conditions:
\begin{itemize}
\item[(a)] $X$ has a uniformly Gâteaux differentiable norm.
\item[(b)] $X$ has Reich's property.
\end{itemize}

\end{remark}


\begin{cor}
Let $H$ be a Hilbert space, $C\subset H$ a nonempty closed convex
subset, $T:C\to C$ a nonexpansive mapping with fixed point set
$F\neq\emptyset$ and $\Phi:C\to C$ a $\rho$-contraction. Then the
approximating curve $\{x_t\}\subset C$ defined by \eqref{Implicit
Iter} strongly converges, as $t\fle0$, to the unique solution $q\in
F$ to the inequality
\begin{equation}\label{inequality Th.-H}
\langle(\Phi-I)q,x-q\rangle\leq0,\,\,\, \forall x\in F.
\end{equation}
\end{cor}

\begin{theorem}
\label{thm:explicit}Let $X$ be a uniformly smooth Banach space, $C$
a nonempty closed convex subset of $X$, $T:C\to C$ a nonexpansive
mapping with $F\neq\emptyset,$ $\Phi:C\to C$ a $\rho$-contraction
and $\{\alpha_{n}\}$ a sequence in $(0,1)$ satisfying
\begin{itemize}
\item[(H1)] $\ds\lim_{n\fle\infty}\alpha_n=0$
\item[(H2)] $\ds\sum_{n=1}^{\infty}\alpha_{n}=\infty$
\item[(H3)] $\ds\sum_{n=1}^{\infty}|\alpha_{n+1}-\alpha_{n}|<\infty$ or
$\ds\lim_{n\fle\infty}\dfrac{\alpha_{n}}{\alpha_{n+1}}=1.$
\end{itemize}
Then, for any $x_0\in C$, the sequence $\{x_{n}\}$ defined by
\begin{equation}\label{Explicit iter}
x_{n+1}=T(\alpha_{n}\Phi(x_{n})+(1-\alpha_{n})x_{n})
\end{equation}
strongly converges to the unique solution $q\in F$ to the inequality
\[ \langle(\Phi-I)q,J(x-q)\rangle\leq0,\,\,\,\forall x\in F.\]
\end{theorem}

\begin{proof}
Since $T$ has fixed points, by Lemma \ref{lemma-useful} (3) we have
that $\{x_n\}$ is bounded, and therefore so are $\{T(x_n)\}$ and
$\{\Phi( x_n)\}$. The fact that $\{x_{n}\}$ is asymptotically
regular is a consequence of Lemmas \ref{lemma-useful} and \ref{lem:xu1}.
Indeed, by hypothesis we have that
$\ds\sum_{n=1}^{\infty}(1-\rho)\alpha_{n-1}=\infty$ and either
$\ds\sum_{n=1}^{\infty}|\alpha_{n}-\alpha_{n-1}|<\infty$
or\begin{equation}
\limsup_{n\fle\infty}\dfrac{|\alpha_{n}-\alpha_{n-1}|}{\alpha_{n}}=
\lim_{n\fle\infty}\left|1-\dfrac{\alpha_{n-1}}{\alpha_{n}}\right|=0.
\end{equation}
Then inequality \eqref{ineq-xn+1-xn}
$$\|x_{n+1}-x_n\| \le(1-(1-\rho)\alpha_n)\|x_n-x_{n-1}\|
 +|\alpha_n-\alpha_{n-1}|\cdot\|\Phi(x_{n-1})-x_{n-1}\|$$
allows us to use Lemma \ref{lem:xu1} to deduce that
\begin{equation}\label{eq:as_reg}
\lim_{n\fle\infty}\|x_{n+1}-x_{n}\|=0.
\end{equation}
Using inequality \eqref{ineq-xn-Txn} and the hypothesis (H1) we get
\begin{equation}\label{eq:x_n-Tx_n}
\lim_{n\fle\infty}\|x_{n}-Tx_{n}\|\leq
\lim_{n\fle\infty}\big(\|x_{n+1}-x_n\| + \alpha_n\|\Phi(x_n)-x_n\|\big)=0.
\end{equation}
We will see now that
\begin{equation}
\limsup_{n}\langle\Phi(q)-q,J(x_{n}-q)\rangle\leq0.\label{eq:limsup}
\end{equation}
Let $\{\beta_{k}\}$ be a null sequence in $(0,1)$ and define
$\{y_{k}\}$ by\[ y_{k}:=T(\beta_{k}\Phi(y_{k})+(1-\beta)y_{k}).\] We
have proved in Theorem \ref{thm:implicit} that $y_{k}$ strongly
converges to $q.$ For any $n,k\geq 0$ define
$$\delta_{n,k}:=\|x_{n}-Tx_{n}\|^{2}+2\|x_{n}-Tx_{n}\|\|y_{k}-Tx_{n}\|$$
and
$$\epsilon_{k}:=\sup_{n}\{\|\Phi(y_{k})-x_{n}\|\|J(\beta_{k}(\Phi(y_{k})-x_{n})+
(1-\beta_k)(y_{k}-x_{n}))-J(y_{k}-x_{n})\|\}.$$
For any fixed $k\in\mathbb{N},$ by \eqref{eq:x_n-Tx_n},
$\lim_{n}\delta_{n,k}=0.$ Moreover $\lim_{k}\epsilon_{k}=0$ because
of the uniform continuity of $J$ over bounded sets. Using inequality
\eqref{subdiff ineqN} we obtain
\[
\begin{array}{ll}
\|y_{k}-x_{n}\|^{2} & \leq(\|Tx_{n}-x_{n}\|+\|y_{k}-Tx_{n}\|)^{2}\\
 & \leq\delta_{n,k}+\|(1-\beta_{k})(y_{k}-x_{n})+\beta_{k}(\Phi(y_{k})-x_{n})\|^{2}\\
 & \leq\delta_{n,k}+(1-\beta_{k})^{2}\|y_{k}-x_{n}\|^{2}+2\beta_{k}
 \langle\Phi(y_{k})-x_{n},J(y_{k}-x_{n})\rangle\\
 & \quad+2\beta_{k}\epsilon_{k}\\
 & =\delta_{n,k}+((1-\beta_{k})^{2}+2\beta_{k})\|y_{k}-x_{n}\|^{2}+2\beta_{k}
 \langle\Phi(y_{k})-y_{k},J(y_{k}-x_{n})\rangle\\
 & \quad+2\beta_{k}\epsilon_{k}.\end{array}\]
Then we deduce that
\begin{equation}
\limsup_{n\fle\infty}\langle\Phi(y_{k})-y_{k},J(x_{n}-y_{k})\rangle\leq
\dfrac{\beta_{k}}{2}\limsup_{n\fle\infty}\|y_{k}-x_{n}\|^{2}+\epsilon_{k}.\label{eq:quasi_dis}
\end{equation}
On the other hand\begin{equation}
\begin{array}{ll}
\langle\Phi(q)-q,J(x_{n}-q)\rangle & =\langle\Phi(q)-q,J(x_{n}-q)-J(x_{n}-y_{k})\rangle\\
 & \quad+\langle(\Phi(q)-q)-(\Phi(y_{k})-y_{k}),J(x_{n}-y_{k})\rangle\\
 & \quad+\langle\Phi(y_{k})-y_{k},J(x_{n}-y_{k})\rangle.\end{array}\label{eq:zero}\end{equation}
Note that \begin{equation}
\lim_{k}(\sup_{n}\{\langle\Phi(q)-q,J(x_{n}-q)-J(x_{n}-y_{k})\rangle\})=0\label{eq:uno}\end{equation}
because $J$ is norm to norm uniform continuous on bounded sets. By
using \eqref{eq:quasi_dis},\eqref{eq:uno} and passing first to
$\limsup_{n}$ and then to $\lim_{k}$, from \ref{eq:zero} we obtain
\[
\limsup_{n}\langle\Phi(q)-q,J(x_{n}-q)\rangle\leq0.\] Finally we
prove that $\{x_{n}\}$ strongly converges to $q.$ Set
$$\eta_{n}:=\|J(\alpha_{n}(\Phi(x_{n})-x_{n})+(x_{n}-q))-J(x_{n}-q)\|.$$
Then $\eta_{n}\to0,$ as $n\fle 0$. We compute\[
\begin{array}{ll}
\|x_{n+1}-q\|^{2} & \leq\|\alpha_{n}(\Phi(x_{n})-q)+(1-\alpha_{n})(x_{n}-q)+\alpha_{n}(\Phi(q)-q)\|^{2}\\
 & \leq\|\alpha_{n}(\Phi(x_{n})-q)+(1-\alpha_{n})(x_{n}-q)\|^{2}\\
 & \quad+2\alpha_{n}\langle\Phi(q)-q,J(x_{n}-q)\rangle\\
 & \quad+2\alpha_{n}\eta_{n}\|\Phi(q)-q\|\\
 & \le(1-(1-\rho)\alpha_{n})\|x_{n}-q\|^{2}+2\alpha_{n}\langle\Phi(q)-q,J(x_{n}-q)\rangle\\
 & \quad+2\alpha_{n}\eta_{n}\|\Phi(q)-q\|\end{array}\]
and the result follows from \eqref{eq:limsup} and Lemma
\ref{lem:xu1}.

\end{proof}

\begin{cor}\label{corollary}
Let $X$ be a uniformly smooth Banach space, $C$ a nonempty closed
convex subset of $X$, $T:C\to C$ a nonexpansive mapping with
$F\neq\emptyset$ and $g:C\to C$ a mapping such that $I-\mu g$ is a
contraction for some $\mu >0$. Assume that $\{\alpha_{n}\}$ is a
sequence in $(0,1)$ satisfying hypotheses (H1)-(H3) in Theorem
\ref{thm:explicit}. Then the sequence $\{x_n\}$ defined by the
iterative scheme
\begin{equation}\label{alg. cor.}
x_{n+1}=Tx_n-\alpha_{n}g(Tx_n),
\end{equation}
strongly converges to the unique solution $p\in F$ to the inequality
problem
\begin{equation}\label{ineq Cor.}
\langle g(p), x-p\rangle \geq 0\,\,\,\forall x\in F.
\end{equation}

\end{cor}
\begin{proof}
Consider the sequence $\{z_n\}$ defined by $z_n=Tx_n$, for any
$n\geq 0$. Then
\begin{eqnarray*}
z_{n+1} & = & T(Tx_n -\alpha_n g(Tx_n))\\
 & = & T(z_n-\frac{\alpha_n}{\mu}\mu g(z_n))\\
 & = & T( \alpha'_n (I-\mu g) z_n+(1-\alpha'_n)z_n),
\end{eqnarray*}
where $\alpha'_n=\frac{\alpha_n}{\mu}$ for all $n\geq 0$, so the
sequence $\{\alpha'_n\}$ satisfies hypotheses (H1)-(H3). Since
$\Phi:=I-\mu g$ is a contraction, Theorem \ref{thm:explicit} implies
the strong convergence of $\{z_n\}$ to the unique solution $p\in F$
to the inequality problem
$$\langle (\Phi-I)p, x-p\rangle \geq 0\,\,\,\forall x\in F,$$
which is equivalent to \eqref{ineq Cor.}. Therefore, from the
iteration scheme \eqref{alg. cor.} we deduce that the sequence
$\{x_n\}$ strongly converges to $p$.
\end{proof}


\section{Rates of asymptotic regularity}\label{s-rate-as-reg}

In the following, we apply proof mining techniques to get effective
rates of asymptotic regularity for the iteration $\{x_n\}$  defined
by (\ref{Exp Iter}). The methods we use in this paper are inspired
by those used in \cite{L07a} to obtain effective rates of asymptotic
regularity for Halpern iterations. As in the case of Halpern
iterations, the main ingredient turns out to be the quantitative
Lemma \ref{quant-Liu-L}.

\begin{theorem}\label{main-thm}
Let $X$ be a normed space, $C\se X$  a nonempty convex subset and
$T:C\to C$ be nonexpansive. Assume that $\Phi:C\to C$ is a
$\rho$-contraction and that $\{\alpha_n\}_{n\geq 0}$ is a sequence
in $(0,1)$  such that $\ds\lim_{n\to\infty} \alpha_n =0$,
$\ds\sum_{n=0}^\infty \alpha_n$ is divergent and
$\ds\sum_{n=0}^\infty|\alpha_{n+1}-\alpha_n|$ is convergent. Let
$x_0\in C$ and $\{x_n\}_{n\ge 0}$ be defined by (\ref{Exp Iter}).
Assume that $\{x_n\}$ is bounded.

Then $\ds\lim_{n\to\infty} \|x_n-Tx_n\|=0$ and moreover for all $\eps\in(0,2)$,\bua
\|x_n-Tx_n\|< \eps,\quad \forall n\ge \Psi(\vp,\beta,\theta,\rho, M,D,\eps), \label{conclusion-main-thm}
\eua where \be
\item $\vp:(0,\infty)\to\Z_+$ is a rate of convergence of $\{\alpha_n\}$,
\item $\beta:(0,\infty)\to\Z_+$ is a Cauchy modulus of $\ds\sum_{n=0}^\infty|\alpha_{n+1}-\alpha_n|$,
\item $\theta:\Z_+\to\Z_+$ is a rate of divergence of $\ds\sum_{n=0}^\infty \alpha_n$,
\item $M\ge 0$ is such that $M\ge \|\Phi(x_0)-x_0\|$,
\item $D>0$ satisfies $D\geq \|x_n-x_m\| \text{~~for all~~} m,n\geq 0$,
\ee and $\Psi(\vp,\beta,\theta,\rho, M,D,\eps)$ is defined by \bua
\ds\Psi:=\max\left\{1+\theta\left(\left\lceil\frac{1}{1-\rho}\right\rceil\left(\beta\left(\frac\eps
{4P}\right)+2+\left\lceil\ln\left(\frac{4D}\eps\right)\right\rceil\right)\right),\,\,
1+\vp\left(\frac\eps {2P}\right)\right\},\label{def-Phi-main-thm}
\eua with $P=(1+\rho)D+M$.
\end{theorem}

\begin{proof}
Applying (\ref{ineq-xn+1-xn}) and (\ref{ineq-Phi-xn}), we get that
for all $n\geq 1$ \bua
\|x_{n+1}-x_n\| &\le & (1-(1-\rho)\alpha_n)\|x_n-x_{n-1}\|\\
&&\,\,+\,|\alpha_n-\alpha_{n-1}|\cdot\|\Phi(x_{n-1})-x_{n-1}\|\\
&\le &
(1-(1-\rho)\alpha_n)\|x_n-x_{n-1}\|+P\cdot|\alpha_n-\alpha_{n-1}|.
\eua Let us denote for $n\geq 1$
\[a_n:=\|x_n-x_{n-1}\|,  \quad b_n:=P\cdot|\alpha_n-\alpha_{n-1}|, \quad \lambda_{n}:=(1-\rho)\alpha_{n-1}.\]
Then $D$ is an upper bound on $\{a_n\}$ and
\[a_{n+1}\leq (1-\lambda_{n+1})a_n+b_n \quad \text{ for all } n\geq 1.\]

Moreover, $\,\ds\sum_{n=1}^\infty \lambda_n$ is divergent with rate
of divergence \beq \delta:\Z_+\to\Z_+, \quad
\delta(n)=1+\theta\left(\left\lceil\frac{1}{1-\rho}\right\rceil\cdot
n\right), \eeq since for all $n\in\Z_+$, \bua
\sum_{i=1}^{\delta(n)}\lambda_i&=&(1-\rho)\sum_{i=0}^{\delta(n)-1}\alpha_i=
(1-\rho)\sum_{i=0}^{\theta\left(\lceil 1/1-\rho\rceil\cdot n\right)}\alpha_i\ge(1-\rho)
\left\lceil\frac{1}{1-\rho}\right\rceil\cdot n\\
&\ge& n. \eua Let $t_n:=\ds\sum_{i=0}^n|\alpha_{i+1}-\alpha_i|$ and
$\ds s_n:=\sum_{i=1}^n b_i=Pt_{n-1}$ and define \beq
\gamma:(0,\infty)\to \Z_+, \quad
\gamma(\eps):=1+\beta\left(\frac\eps{P}\right), \eeq Then for all
$n\ge 0$, \bua \ds s_{\gamma(\eps)+n}-s_{\gamma(\eps)}&= &
P\left(t_{\beta\left(\eps/P\right)+n}-t_{\beta\left(\eps/P\right)}\right)<
P\cdot\frac\eps{P}=\eps. \eua Thus, $\ds \sum_{n=1}^\infty b_n$ is
convergent with Cauchy modulus $\gamma$.

It follows that we can apply Lemma \ref{quant-Liu-L} to get that for
all $\eps\in(0,2)$ and for all $n\geq
h_1(\beta,\theta,\rho,M,D,\eps)$ \beq \|x_n-x_{n-1}\|<\frac\eps
2,\label{final-1} \eeq where \bua
h_1(\beta,\theta,\rho,M,D,\eps)&:=& 1+\theta\left(\left\lceil\frac{1}{1-\rho}
\right\rceil\cdot\left(\beta\left(\frac\eps {4P}\right)+2+\left\lceil\ln
\left(\frac{4D}\eps\right)\right\rceil\right)\right)\\
\eua Using (\ref{ineq-xn-Txn}) and (\ref{ineq-Phi-xn}), we get that
for all $n\geq 1$, \beq \ba{lll}
\|x_{n-1}-Tx_{n-1}\| & \leq & \|x_n-x_{n-1}\|+\alpha_{n-1}\|\Phi(x_{n-1})-x_{n-1}\|\\
&\le & \|x_n-x_{n-1}\| + P\alpha_{n-1}. \ea\label{final-0} \eeq Let
$\ds h_2(\vp,\rho,M,D,\eps):=1+\vp\left(\frac\eps {2P}\right)$.
Since $\vp$ is a rate of convergence of $\{\alpha_n\}$ towards $0$, it
follows that \beq P\alpha_{n-1}< \frac\eps 2 \quad \text{ for all }
n\geq h_2(\vp,\rho,M,D,\eps).\label{final-2} \eeq

As a consequence of (\ref{final-1}), (\ref{final-0}) and
(\ref{final-2}), we get that
\[\|x_{n-1}-Tx_{n-1}\|<\eps\]
for all $n\geq \max\{h_1(\beta,\theta,\rho,M,D,\eps),
h_2(\vp,\rho,M,D,\eps)\}$, so the conclusion of the theorem follows.
\end{proof}

If $C$ is bounded, then $\{x_n\}$ is bounded for all $x_0\in C$.
Moreover, we can take $M:=D:=d_C$ in the above theorem, where
$d_C:=\sup\{\|x-y\|\mid x,y\in C\}$ is the diameter of $C$.

\begin{cor}\label{C-bounded}
In the hypotheses of Theorem \ref{main-thm}, assume moreover that
$C$ is bounded.

Then  $\ds\lim_{n\to\infty} \|x_n-Tx_n\|=0$ for all $x_0\in C$ and
moreover for all $\eps\in(0,2)$,
\bua
\|x_n-Tx_n\|< \eps,\quad \forall n\ge \Psi(\vp,\beta,\theta,\rho,d_C,\eps),
\eua
where $\Psi(\vp,\beta,\theta,\rho,d_C,\eps)$ is defined as in
Theorem \ref{main-thm} by replacing $M$ and $D$ with $d_C$.
\end{cor}

Thus, for bounded $C$, we get asymptotic regularity for general
$\{\alpha_n\}$ and an explicit rate of asymptotic regularity
$\Psi(\vp,\beta,\theta,\rho,d_C,\eps)$  that depends weakly on $C$
(via its diameter $d_C$) and on the $\rho$-contraction $\Phi$ (via
$\rho$), while it does not depend on the nonexpansive mapping $T$,
the starting point $x_0\in C$ of the iteration or other data related
with $C$ and $X$.

The rate of asymptotic regularity can be simplified when the
sequence $\{\alpha_n\}$ is decreasing.

\begin{cor}\label{lambda-decreasing}
Let $X,C,T,\Phi,\{x_n\}$ be as in the hypotheses of Corollary
\ref{C-bounded}. Assume that $\{\alpha_n\}$ is a decreasing sequence
in $(0,1)$  such that $\ds\lim_{n\to\infty} \alpha_n =0$ and
$\ds\sum_{n=0}^\infty \alpha_n$ is divergent.

Then  $\ds\lim_{n\to\infty}  \|x_n-Tx_n\|=0$ for all $x_0\in C$ and
furthermore, for all  $ \eps\in(0,2)$,
\bua
\|x_n-Tx_n\|< \eps, \quad \forall n\ge \Psi(\vp,\theta,\rho,d_C,\eps),
\eua where $\vp:(0,\infty)\to\Z_+$ is a rate of convergence of
$\{\alpha_n\}$, $\theta:\Z_+\to\Z_+$ is a rate of divergence of
$\,\ds\sum_{n=0}^\infty \alpha_n$ and $\Psi(\vp,\theta,\rho,
d_C,\eps)$ is defined by \bua
\begin{array}{l}
\ds\Psi:=\max\left\{1+\theta\left(\left\lceil\frac{1}{1-\rho}\right\rceil\left(\vp\left(\frac\eps
{4P}\right)+2+\left\lceil\ln\left(\frac{4d_C}\eps\right)\right\rceil\right)\right),\,\,
1+\vp\left(\frac\eps {2P}\right)\right\}
\end{array}
\eua with $P=(2+\rho)d_C$.
\end{cor}

\begin{proof}
Remark that $\ds\sum_{n=0}^\infty|\alpha_{n+1}-\alpha_n|$ is
convergent with Cauchy modulus $\vp$.
\end{proof}

Finally, we get, as in the case of Halpern iterates,  an exponential
(in $1/\eps$) rate of asymptotic regularity for $\ds\alpha_n=\frac
1{n+1}$.

\begin{cor}\label{lambda-1-n+1}
Let $X,C,T,\Phi,\{x_n\}$ be as in the hypotheses of Corollary
\ref{C-bounded}. Assume that $\ds\alpha_n=\frac 1{n+1}$ for all
$n\geq 0$.

Then  $\ds\lim_{n\to\infty} \|x_n-Tx_n\|=0$ for all $x_0\in C$ and
furthermore, for all  $ \eps\in(0,2)$,   \bua
\|x_n-Tx_n\|< \eps,\quad\forall n\ge \Theta(\rho,d_C,\eps),
\eua where \bua
\begin{array}{l}
\Theta(\rho,d_C,\eps)=\exp\left(\ds\frac{4}{1-\rho}\left(\frac{16d_C}{\eps}+2\right)\right)
\end{array}
\eua
\end{cor}

\begin{proof}
We can apply Corollary \ref{lambda-decreasing} with
\[\vp:(0,\infty)\to\Z_+, \quad \vp(\eps)= \left\lceil\frac{1}\eps\right\rceil-1.\]
and
\[\theta:\Z_+\to\Z_+, \quad \theta(n)=4^{n}-1\]
to conclude that for all $\eps\in(0,2)$,
\[\|x_n-Tx_n\|< \eps,\quad \forall n\ge \Psi(\rho,d_C,\eps),\]
where $P=(2+\rho)d_C$ and \bua
\ds\Psi(\rho,d_C,\eps) &=&\max\left\{\exp\left(\ln 4\cdot \left\lceil
\frac{1}{1-\rho}\right\rceil\cdot\left(\left\lceil\frac{4P}{\eps}
\right\rceil+1+\left\lceil\ln\left(\frac{4d_C}\eps\right)\right\rceil\right)\right)
,\,\, \left\lceil\frac{2P}\eps \right\rceil\right\}\\
&=& \exp\left(\ln 4\cdot \left\lceil\frac{1}{1-\rho}\right\rceil\cdot\left
(\left\lceil\frac{4P}{\eps} \right\rceil+1+\left\lceil\ln\left(\frac{4d_C}\eps
\right)\right\rceil\right)\right)\\
& < &
\exp\left(\frac{4}{1-\rho}\left(\frac{16d_C}{\eps}+2\right)\right) =
\Theta(\rho,d_C,\eps). \eua as $\rho\in (0,1)$, $\lceil a\rceil <
a+1$ and $1+\ln a\leq a$ for all $a>0$. The conclusion follows now
immediately.
\end{proof}

\section{Applications}
As it was pointed out in the introduction, iterative methods for
nonexpansive mappings have been applied to solve the problem of
finding a solution to the VIP($A,C$) \eqref{VIP-intro} which, in
fact, is equivalent under suitable conditions to finding the
minimizer of a certain function. On the other hand, the relation
between the set of zeros of an accretive operator and the fixed
point set of its resolvent allows us to use those iterative
techniques for nonexpansive mappings to approximate zeros of such
operators. We first apply the explicit iterative method for
approximating fixed points, presented in section 3, to solve a
particular variational inequality problem in the setting of Hilbert
spaces. Then we focus on the asymptotic behavior of the resolvent of
an accretive operator in the framework of Banach spaces.

\subsection{A variational inequality problem}
Let $H$ be a Hilbert space, 
$T:H\to H$ be a nonexpansive mapping with fixed point set
$F\neq\emptyset$, and $\Phi:H\fle H$ be a contraction. Assume that
$A$ is a Lipschitzian self-operator on $H$ which is strongly
monotone; that is, there exist a constant $\eta >0$ such that
$$\langle Ax-Ay ,x-y\rangle\geq \eta \Vert x-y\Vert^2, \,\,\forall x,y\in H.$$
It is known that the following variational inequality
\begin{equation}\label{VIP-appl}
\langle(A-\gamma \Phi) q , q-x\rangle \leq 0\,\,\, \forall x\in F,
\end{equation}
where $\gamma>0$, is the optimality condition for the minimization
problem
$$ \min_{x\in F} f(x)- h(x)$$
where $f$ is a differentiable function with differential $\partial
f=A$ and $h$ is a potential function for $\gamma \Phi$ (i.e.
$h'(x)=\gamma \Phi(x)$ for $x\in H$). Marino and Xu \cite{MX06}
presented an iterative method to solve the variational inequality
\eqref{VIP-appl} for a linear bounded operator. We apply the
explicit method \eqref{Explicit iter}, in particular algorithm
\eqref{alg. cor.}, to solve such variational inequality dispensing
with the linear condition on the operator $A$. For that purpose,
we need the following lemma.

\begin{lemma}\label{Lemma-Appl}
Assume that $A$ is a $L$-Lipschitzian $\eta$-strongly monotone
operator, and let $\Phi$ be a $\rho$-contraction. Then, for any
$\gamma<\eta/\rho$, $A-\gamma \Phi$ is $R$-Lipschitzian and
$\delta$-strongly monotone with $R=L+\gamma\rho$ and
$\delta=\eta-\gamma\rho$. Besides, for any $0<\mu<2\delta/R^2$, the
mapping $I- \mu( A-\gamma\Phi)$ is a contraction.
\end{lemma}

\begin{proof}
Since $A$ is $L$-Lipschitzian and $\Phi$ is a $\rho$-contraction,
$$\|(A-\gamma\Phi)x-(A-\gamma\Phi)y\|\leq\|Ax-Ay\|+\gamma\|\Phi x-\Phi y\|\leq (L+\gamma\rho)\|x-y\|,$$
that is, $A-\gamma \Phi$ is Lipschitzian with constant
$R=L+\gamma\rho$. The strong monotonicity of $A-\gamma \Phi$ is
consequence of the strong monotonicity of $A$ as it is showed as
follow.
\begin{eqnarray*}
\langle (A-\gamma\Phi)x-(A-\gamma\Phi)y,x-y\rangle & = & \langle Ax-Ay,x-y\rangle
- \gamma \langle \Phi x-\Phi y,x-y\rangle\\
 & \geq & \eta\|x-y\|^2-\gamma\|\Phi x-\Phi y\| \|x-y\|\\
  & \geq & (\eta-\gamma\rho) \|x-y\|^2,
\end{eqnarray*}
where $\delta=\eta-\gamma\rho>0$. By applying the $R$-Lipschitz
continuity and $\delta$-strong monotonicity of $B:=A-\gamma \Phi$ we
obtain
\begin{eqnarray*}
\|(I-\mu B)x-(I-\mu B)y\|^2 & = &
\|x-y\|^2+\mu^2\|Bx-By\|^2-\mu\langle x-y,Bx-By\rangle\\
 & \leq & \|x-y\|^2+\mu^2R^2\|x-y\|^2-2\mu\delta \|x-y\|^2\\
 & = & (1-\mu(2\delta-\mu R^2))\|x-y\|^2.\\
\end{eqnarray*}
Then, for any $0<\mu<2\delta/R^2$, the mapping $I-\mu
(A-\gamma\Phi)$ is a contraction with constant
$\sqrt{1-\mu(2\delta-\mu R^2)}$.
\end{proof}

\begin{theorem} Let $T$ be a nonexpansive mapping with fixed point set $F$,
$A$ a $L$-Lipschitzian $\eta$-strongly monotone operator and $\Phi$
a $\rho$-contraction on a Hilbert space. Then, for any
$\gamma<\eta/\rho$, the sequence defined by the iterative scheme
$$ x_{n+1} =  Tx_n -\alpha_{n} (A-\gamma \Phi) Tx_n,$$
where $\{\alpha_n\}\subset(0,1)$ satisfies hypotheses (H1)-(H2) in
Theorem \ref{thm:explicit}, strongly converges to the unique
solution to the variational inequality \eqref{VIP-appl}.
\end{theorem}

\begin{proof}
Note that, for any $\gamma<\eta/\rho$, Lemma \ref{Lemma-Appl}
implies that there exists $\mu>0$ such that $I-\mu g$ is a
contraction, where $g=A-\gamma\Phi$. Then, by Theorem
\ref{corollary} we obtain the strong convergence of the sequence
$\{x_n\}$ to the unique solution to the variational inequality
problem \eqref{VIP-appl}.

\end{proof}

\subsection{Zeros of $m$-accretive operators}

Let $X$ be a real Banach space. A set-valued operator $A:
X\rightarrow 2^X$ with domain $D(A)$ and range $R(A)$ in $X$ is said
to be \emph{accretive} if, for each $x_i\in D(A)$ and $y_i\in Ax_i$
($i=1,2$), there exists $j(x_1-x_2)\in J(x_1-x_2)$ such that
$$\langle y_1-y_2,j(x_1-x_2)\rangle\geq0,$$
where $J$ is the normalized duality mapping. An accretive operator
$A$ is \emph{m-accretive} if $R(I+\lambda A)=X$ for all
$\lambda>0$. Denote the set of zeros of $A$ by
$$Z:=A^{-1}(0)=\{z\in D(A):0\in Az\}.$$

Throughout this subsection it is assumed that $A$ is $m$-accretive
and $A^{-1}(0)\neq\emptyset$. Set $C=\overline{D(A)}$ and assume it
is convex. It is known that the \emph{resolvent} of $A$, defined by
$$J_\lambda=(I+\lambda A)^{-1},$$
for $\lambda>0$, is a single-valued nonexpansive mapping from $C$
into itself (cf. \cite{B76}).

If we consider the problem of finding a zero of $A$, i.e.,
$$\textrm{ find } z\in C \textrm{ such that } 0\in Az,$$
it is straightforward to see that the set of zeros $A^{-1}(0)$
coincides with the fixed point set of $J_\lambda$,
$Fix(J_\lambda)$, for any $\lambda>0$. Therefore an equivalent
problem is to find $z\in Fix(J_\lambda)$.

As a consequence of the convergence of the implicit iterative scheme
\eqref{Implicit Iter}, we obtain Reich's result (cf. \cite{R80}) for
approximating zeros of accretive operators in uniformly smooth
Banach spaces. Besides, the following theorem in the setting of
reflexive Banach spaces with weakly continuous normalized duality
mapping constitutes a new approach.

\begin{theorem} Let $A$ be a m-accretive operator in either a
uniformly smooth Banach space or a reflexive Banach spaces with
weakly continuous normalized duality mapping $X$. Then, for each
$x\in X$, the sequence $\{J_\lambda (x)\}$ strongly converges, as
$\lambda\fle\infty$, to the unique zero of $A$, $q\in A^{-1}(0)$,
which satisfies the variational inequality
\begin{equation}\label{sunny nonexp. ineq}
\langle x- q,J( y-q)\rangle\leq0\,\,\,\forall y\in A^{-1}(0).
\end{equation}
\end{theorem}
\begin{proof}
Given $x\in X$ we consider the approximating curve $\{x_t\}$ such
that $x_t=J_{1/t}x$, for any $t\in (0,1)$. By definition of the
resolvent of $A$, we obtain the following equivalence:
\begin{eqnarray*}
x_t=(I+\frac{1}{t}A)^{-1}x & \Leftrightarrow & x\in
x_t+\frac{1}{t}Ax_t\\
 & \Leftrightarrow & t(x-x_t)\in Ax_t\\
 & \Leftrightarrow & x_t +t(x-x_t)\in (I+A)x_t\\
 & \Leftrightarrow & x_t = (I+A)^{-1}(x_t+t(x-x_t))\\
 & \Leftrightarrow & x_t= T(t \Phi(x_t)+(1-t)x_t),
\end{eqnarray*}
where $T=(I+A)^{-1}$ is the nonexpansive resolvent with constant
$1$, and $\Phi=x$ is a constant mapping which is a contraction.
Therefore, Theorem \ref{thm:implicit} implies the strong convergence
of $\{x_t\}$, as $t\fle 0$, to the unique solution to the inequality
\eqref{inequality Th.}; in other words, $\{J_\lambda x\}$ strongly
converges, as $\lambda\fle\infty$, to the unique solution $q\in
A^{-1}(0)$ to the inequality \eqref{sunny nonexp. ineq}.
\end{proof}

\begin{remark}
If we define the mapping $Q:X\fle A^{-1}(0)$ such that, for any
$x\in X$,
$$Qx=\lim_{\lambda\to\infty}J_{\lambda}x,$$
then, since $Qx$ satisfies the inequality \eqref{sunny nonexp.
ineq}, by Lemma \ref{sunny} we can claim that $Q$ is the unique
sunny nonexpansive retraction from $X$ to $A^{-1}(0)$.

\end{remark}

\end{document}